\RequirePackage[l2tabu]{nag}
\documentclass[11pt]{article}

\usepackage[T1]{fontenc}
\usepackage[utf8]{inputenc}
\usepackage{lmodern}

\usepackage{amssymb,amsthm}
\usepackage{mathtools}

\usepackage[hidelinks]{hyperref}
\usepackage[vmargin=2cm,hmargin=3cm,includefoot]{geometry}

\usepackage{xcolor}
\hypersetup{
    colorlinks,
    linkcolor={red!50!black},
    citecolor={blue!50!black},
    urlcolor={blue!80!black}
}

\newcommand*{\N}{\mathbb{N}}
\newcommand*{\Z}{\mathbb{Z}}
\newcommand*{\R}{\mathbb{R}}

\newcommand*{\eps}{\varepsilon}

\newcommand*{\vol}{\mathrm{vol}}

\newcommand*{\calL}{\mathcal{L}}
\newcommand*{\calC}{\mathcal{C}}
\newcommand*{\calZ}{\mathcal{Z}}

\newcommand*{\frakX}{\mathfrak{X}}

\DeclarePairedDelimiter\abs{\lvert}{\rvert}
\DeclarePairedDelimiter\norm{\lVert}{\rVert}
\DeclarePairedDelimiterX\innerp[2]{\langle}{\rangle}{#1,#2}
\DeclarePairedDelimiterX\intervCO[2]{[}{[}{#1,#2}

\DeclareMathOperator{\id}{id}
\let\div\relax\DeclareMathOperator{\div}{div}
\DeclareMathOperator{\grad}{\nabla}
\DeclareMathOperator{\Hess}{Hess}
\DeclareMathOperator{\lapla}{\bigtriangleup}
\DeclareMathOperator{\levici}{\nabla^{\textsc{lc}}}

\DeclareMathOperator{\extdiff}{d}
\DeclareMathOperator{\Lie}{\calL}

\DeclareMathOperator{\interior}{int}
\DeclareMathOperator{\Ric}{Ric}

\newtheorem*{theorem*}{Theorem}
\newtheorem*{divthm}{Divergence theorem}
\theoremstyle{definition}

\theoremstyle{remark}

\newtheorem*{remarks*}{Remarks}

\author{Beno\^{\i}t Jubin}
\title{Closed Kac--Rice type formulas on Riemannian manifolds}
\date{\today}

\begin{document}

\maketitle

\begin{abstract}
We establish a few formulas that compute the volume of the zero-set (or nodal set) of a function on a compact Riemannian manifold as integrals of functionals of the function and its derivatives.
\end{abstract}

\section{Introduction}
In their study of the volume of the zero-sets (or nodal sets) of Gaussian random fields on $(\R/\Z)^n$, the authors of~\cite{nodal} established a few formulas (of ``Kac--Rice type'') that compute the volume of the zero-set of a real-valued function $f$ on $(\R/\Z)^n$ as an integral of a functional of $f$ and its derivatives.
In particular, they gave a general formula in the one-dimensional case and a few specific formulas in higher dimension.
In this note, we establish a general formula for functions on compact Riemannian manifolds.

The importance of these formulas for the applications studied in~\cite{nodal}, compared to existing Kac--Rice formulas, stems from the fact that they are in ``closed form'' as opposed to being limits of an integral depending on a parameter.
Also, in order to apply techniques of the Malliavin calculus to these formulas, one needs the integrands to be Lipschitz continuous functionals of $f$ and its derivatives, a requirement met by Formula~\eqref{eq:lipschitz}.

The only difference between~\eqref{eq:lipschitz} and~\cite[Prop.~7]{nodal} is the additional term involving the Ricci curvature, $\abs f \eta_f^{-3} \Ric(\grad f, \grad f)$, and this term is in the required domain of the Malliavin calculus by the same proof as \cite[Lem.~2 p.~26]{nodal}.
Therefore, \cite[Thm.~1]{nodal} holds on any compact Riemannian manifold.
Similarly, Formula~\eqref{eq:corner} shows that the extra boundary terms are not problematic, so~\cite[Thm.~1]{nodal} holds on any compact Riemannian manifold with corners, a generalization which includes~\cite[Thm.~2]{nodal} as a special case.

\paragraph{Acknowledgments}
I would like to thank Guillaume Poly for fruitful discussions.

\section{Statement of the results}

Let $(M, g)$ be a compact Riemannian manifold.
If $f \in \calC^2(M,\R)$, we denote by $\sigma_f \colon M \to \{-1, 0, +1 \}$ its sign and we set
\begin{equation}
\eta_f \coloneqq \sqrt{f^2 + \norm{\grad f}^2}.
\end{equation}
If $f$ is nondegenerate, that is, $f(x) = 0$ implies $df(x) \neq 0$, then $\eta_f \in \calC^1(M, \R_{>0})$.
We denote the zero-set of $f$ by
\begin{equation}
\calZ_f \coloneqq f^{-1}(0).
\end{equation}
If~$f$ is nondegenerate, then $\calZ_f$ is a compact Riemannian submanifold of~$M$ of class~$\calC^2$ and codimension~1.
The remaining notation should be clear and will be further explained in the next section.
Our results can be summarized as follows.

\begin{theorem*}
Let $(M, g)$ be a compact $n$-dimensional Riemannian manifold.
Let $f \in \calC^2(M, \R)$ be nondegenerate.
\begin{enumerate}
\item
Let $F \in \calC^1(TM, TM)$ be such that $F(u) \sim_\infty \frac{u}{\norm{u}}$ and $\div \left( F \circ \frac{\grad f}{f} \right) \in \calL^1(M)$.
Then,
\begin{equation}
\label{eq:div} 
\vol(\calZ_f) = -\frac12 \int_M \left( \div \left( F \circ \frac{\grad f}{f} \right) \right) \vol_M.
\end{equation}
\item
Let $G \in \calC^1(\R_{\geq 0}, \R)$ be such that $G(x) \sim_{+\infty} x^{-1}$ and 
$\div \left( \left( G \circ \frac{\norm{\grad f}}{\abs f} \right) \frac{\grad f}{f} \right) \in \calL^1(M)$.
Then,
\begin{multline}
\label{eq:general} 
\vol(\calZ_f) = \frac12 \int_M \left(
\left( G \circ \frac{\norm{\grad f}}{\abs f} \right) \left( \frac{\norm{\grad f}^2}{f^2} - \frac{\lapla f}{f} \right) + {}
\right.\\\left.
\sigma_f \left( G' \circ \frac{\norm{\grad f}}{\abs f} \right) \left( \frac{\norm{\grad f}^3}{f^3} - \frac{\Hess(f)(\grad f, \grad f)}{f^2 \norm{\grad f}} \right)
\right) \vol_M.
\end{multline}

\item
Let $g \in \calC^1(\R_{\geq 0}, \R)$ be such that $\lim_{x \to +\infty} g(x) = 1$ and $g' \in \calL^1(\R)$.
Then,
\begin{multline}
\label{eq:1}
\vol(\calZ_f) = \frac12 \int_M \left(
\frac{\sigma_f}{\eta_f^3} \left( g \circ \frac{\norm{\grad f}}{\abs f} \right)
\left( f \norm{\grad f}^2 + \Hess(f)(\grad f, \grad f) - \eta_f^2 \lapla f \right) + {}
\right.\\\left.
\frac{\norm{\grad f}}{f^2 \eta_f} \left( g' \circ \frac{\norm{\grad f}}{\abs f} \right)
\left( \norm{\grad f}^2 - \frac{f \Hess(f)(\grad f, \grad f)}{\norm{\grad f}^2} \right)
\right) \vol_M.
\end{multline}

Let $g \in \calC^1(\R_{\geq 0}, \R)$ be such that $\lim_{x \to +\infty} g(x) = 1$ and $g' \in \calL^1(\R)$ and $g$ is twice differentiable at 0.
Then,
\begin{multline}
\label{eq:2}
\vol(\calZ_f) = \frac12 \int_M \left(
\frac{\sigma_f}{\norm{\grad f}} \left( g \circ \frac{\norm{\grad f}}{\abs f} \right)
\left( \frac{\Hess(f)(\grad f, \grad f)}{\norm{\grad f}^2} - \lapla f \right) + {}
\right.\\\left.
\left( g' \circ \frac{\norm{\grad f}}{\abs f} \right)
\left(  \frac{\norm{\grad f}^2}{f^2} - \frac{\Hess(f)(\grad f, \grad f)}{f \norm{\grad f}^2} \right)
\right)
\vol_M.
\end{multline}

\item
One has
\begin{equation}
\label{eq:algebraic}
\vol(\calZ_f) = \frac12 \int_M
\frac{\sigma_f}{\eta_f^3} \left(
f \norm{\grad f}^2 + \Hess(f)(\grad f, \grad f)
- \eta_f^2 \lapla f \right)
\vol_M
\end{equation}
and
\begin{multline}
\label{eq:arctan}
\vol(\calZ_f) = \frac1\pi \int_M  \left( 
\norm{\grad f}^{-1} \left( \arctan \circ \frac{\norm{\grad f}}{f} \right)
\left( \frac{\Hess(f)(\grad f, \grad f)}{\norm{\grad f}^2} - \lapla f \right) + {}
\right.\\\left.
\eta_f^{-2} \left( \norm{\grad f}^2 - \frac{f \Hess(f)(\grad f, \grad f)}{\norm{\grad f}^2} \right)
\right)
\vol_M
\end{multline}
and
\begin{multline}
\label{eq:tanh}
\vol(\calZ_f) = \frac12 \int_M \left(
\norm{\grad f}^{-1} \left( \tanh \circ \frac{\norm{\grad f}}{f} \right)
\left( \frac{\Hess(f)(\grad f, \grad f)}{\norm{\grad f}^2} - \lapla f \right) + {}
\right.\\\left.
\left( \cosh \circ \frac{\norm{\grad f}}{f} \right)^{-2}
\left( \frac{\norm{\grad f}^2}{f^2} - \frac{\Hess(f)(\grad f, \grad f)}{f \norm{\grad f}^2} \right)
\right)
\vol_M.
\end{multline}
\end{enumerate}
\end{theorem*}

\begin{remarks*}\quad

\begin{enumerate}
\item
By ``$F(u) \sim_\infty \frac{u}{\norm{u}}$'', we mean that $\lim_{\norm{u} \to +\infty} d \left( F(u), \frac{u}{\norm{u}} \right) =  0$, where $d$ is the distance on $TM$ induced by the Riemannian metric of $M$ (or any distance, since $M$ is compact and $\frac{u}{\norm u}$ has unit norm).

\item
Since $\calZ_f$ is negligible and $\div \left( F \circ \frac{\grad f}{f} \right)$ is defined on $M \setminus \calZ_f$, it makes sense to require ``$\div \left( F \circ \frac{\grad f}{f} \right) \in \calL^1(M)$'' and similarly for Item~2.

\item
In the last three formulas, all terms of the integrands are bounded on $M$ and continuous on $M \setminus \calZ_f$.
Indeed, the Hessian expressions are quadratic in $\norm{\grad f}$, the $\arctan$ and $\tanh$ expressions are linear in $\norm{\grad f}$ when $\norm{\grad f}$ is small, and the $\cosh$ expression is exponentially small in $\abs{f}$ when $\abs{f}$ is small.
However, not all terms need be continuous on $M$.
This problem is treated in the next section.

\item
The cases considered in~\cite{nodal} correspond to $M = (\R/\Z)^n$ with the standard flat metric.
In particular, Formula~\eqref{eq:algebraic} is essentially~\cite[Prop.~5]{nodal} (in the case of $(\R/\Z)^n$ with the standard flat metric).
 
\item
In dimension~1, the general formula~\eqref{eq:div} reduces to~\cite[Prop.~2]{nodal}, and in that case, only the condition $\lim_{x \to \pm\infty} F(x) = \pm 1$ is required if one considers the integral as an improper Lebesgue integral.
Similarly, Formula~\eqref{eq:algebraic} reduces to~\cite[Prop.~1]{nodal} and Formula~\eqref{eq:arctan} to~\cite[Cor.~1 of Prop.~2]{nodal}.

\item
It is possible to extend these results to compact Riemannian manifolds with corners, assuming that $\calZ_f$ intersects $\partial M$ transversely.
The proof goes along the same lines as below, using a divergence theorem on compact manifolds with corners.
Boundary terms will appear in the formulas. For instance, Formula~\eqref{eq:div} becomes
\begin{equation}
\label{eq:corner}
\vol(\calZ_f) = \frac12 \left( \int_{\partial M} \innerp*{F \circ \frac{\grad f}{f}}{v} \vol_{\partial M} - \int_M \left( \div \left( F \circ \frac{\grad f}{f} \right) \right) \vol_M \right).
\end{equation}
This formula reduces in dimension~1 to~\cite[Prop.~3]{nodal}.

\item
Some degenerate functions are considered in~\cite[Prop.~4]{nodal}.
For instance, the 1-dimensional version of Formula~\eqref{eq:arctan} holds for any $f \in \calC^2(\R/\Z,\R)$ such that for any $x \in \R/\Z$ there exists $r \in \N$ such that $f^{(r)}(x)$ exists and is nonzero.
\end{enumerate}
\end{remarks*}

\section{Proof of the results}

Let $(M, g)$ be a compact $n$-dimensional Riemannian manifold with boundary.
Its metric will also be denoted by $\innerp{-}{-}$ and the associated norm by $\norm{-}$.
We denote by $\levici$ its Levi-Civita connection.
Let $\vol_M$ be the Riemannian density on $M$ and $\vol_{\partial M}$ be the induced density on $\partial M$.
Let $v \in \Gamma(TM|_{\partial M})$ be the unit outward normal vectorfield.
The symbol $\vol$ will denote the volume of $(n-1)$-dimensional submanifolds.

The gradient of a function $f \in \calC^1(M, \R)$ is defined by $\grad f \coloneqq (\extdiff f)^\sharp$.
The divergence of a $\calC^1$-vectorfield $X \in \frakX(M)$ is defined by $\Lie_X \vol_M = (\div X) \vol_M$ where $\Lie$ denotes the Lie derivative.
The Laplacian of $\calC^2$-functions is defined by $\lapla \coloneqq \div \circ \grad$.
The Hessian of a function $f \in \calC^2(M, \R)$ is defined by $\Hess f \coloneqq \levici \extdiff f = \innerp{\levici \grad f}{-}$.

Our basic tool is:

\begin{divthm}
Let $(M, g)$ be a compact $n$-dimensional Riemannian manifold with boundary.
Let $X \in \frakX(M)$ be a continuous vectorfield on $M$ which is of class $\calC^1$ on $\interior{M}$ and such that $\div X \in \calL^1(\interior{M})$.
Then,
\begin{equation}
\int_M (\div X) \vol_M = \int_{\partial M} \innerp*{X}{v} \vol_{\partial M}.
\end{equation}
\end{divthm}

\begin{proof}
If $X$ is of class $\calC^1$ on $M$, then this is the standard divergence theorem.
Else, we consider the geodesic flow from the boundary of $M$ along the unit normal vectorfield.
For $\eps > 0$ small enough, set $\theta_\eps \colon \partial M \to M, x \mapsto \exp(x,-\eps v_x)$ and set $M_\eps \coloneqq M \setminus \bigcup_{s \in \intervCO{0}{\eps}} \theta_s(\partial M)$.
For $\eps$ small enough, $M_\epsilon$ is a compact submanifold with boundary of $M$, and $\theta_\eps$ is a diffeomorphism onto $\partial M_\eps$.
Applying the standard divergence theorem on $M_\eps$, one obtains
$\int_{M_\eps} (\div X) \vol_M = \int_{\partial M_\eps} \innerp*{X}{v} \vol_{\partial M_\eps}$.
When $\eps \to 0$, the left-hand side converges to $\int_{M_\eps} (\div X) \vol_M$ by Lebesgue's dominated convergence theorem, since $\div X \in \calL^1(\interior{M})$.
The right-hand side is equal, by change of variable, to $\int_{\partial M} \innerp*{\theta_\eps^* X}{v} (\det T\theta_\eps) \vol_{\partial M}$, which converges to $\int_{\partial M} \innerp*{X}{v} \vol_{\partial M}$ since the integrand is uniformly convergent and $\partial M$ is compact.
\end{proof}

We can now prove the main theorem.

\begin{proof}[Proof of the main theorem]
Let $(M, g)$ be a compact $n$-dimensional Riemannian manifold (without boundary).
Let $f \in \calC^2(M, \R)$ be a nondegenerate function on $M$.
Its zero-set $\calZ_f$ is a compact Riemannian submanifold of dimension $n-1$ of class $\calC^2$.
It has finitely many connected components, say $Z_1, \ldots, Z_d$.
Let $(A_j)_{j \in J}$ be the connected components of $M \setminus \calZ_f$.
They are open submanifolds of $M$ and there are at most $|\pi_0(M)| + 2d$ of them.

\bigskip
\noindent
\textbf{Proof of 1.}
Let $F \in \calC^1(TM, TM)$ be as in the statement.
The vectorfield $F \circ \frac{\grad f}{f}$ is defined on $M \setminus \calZ_f$.
For all $j \in J$, its restriction to $A_j$ extends continuously to $\overline{A_j} = A_j \cup \bigcup_{k \in b(j)} Z_k$ and is equal to the unit inward normal vectorfield on the boundary.
The divergence theorem applied to this extension on $\overline{A_j}$ gives
\begin{align*}
\int_{A_j} \div \left( F \circ \frac{\grad f}{f} \right) \vol_M
&= \int_{\partial A_j} \innerp*{F \circ \frac{\grad f}{f}}{v} \vol_{\partial A_j}\\
&= \int_{\partial A_j} (-1) \vol_{\partial A_j}\\
&= - \sum_{k \in b(j)} \vol(Z_k).
\end{align*}
Therefore, summing over $j \in J$, one obtains
\begin{equation}
\int_M \div \left( F \circ \frac{\grad f}{f} \right) \vol_M = -2 \vol(\calZ_f)
\end{equation}
the factor 2 coming from the fact that each $Z_i$ borders exactly two $A_j$'s.

\bigskip
\noindent
\textbf{Proof of 2.}
We apply (1.) to a function $F$ that is radial.
Namely, let $G$ be as in the statement and apply (1.) to $F \coloneqq \left( G \circ \norm{-} \right) \id_{TM}$.
The Leibniz rule for the divergence gives
\begin{align*}
\div \left( F \circ \frac{\grad f}{f} \right)
&= \div \left( \left(G \circ \frac{\norm{\grad f}}{\abs f} \right) \frac{\grad f}{f} \right)\\
&= \left( G \circ \frac{\norm{\grad f}}{\abs f} \right) \frac{\div \grad f}{f} + \innerp*{\grad \left( \frac1f G \circ \frac{\norm{\grad f}}{\abs f} \right)}{\grad f}\\
&\!\begin{multlined}
= \left( G \circ \frac{\norm{\grad f}}{\abs f} \right) \left( \frac{\lapla f}{f} - \frac{\norm{\grad f}^2}{f^2} \right) +{}\\
\frac1f \left( G' \circ \frac{\norm{\grad f}}{\abs f} \right) \innerp*{\grad \frac{\norm{\grad f}}{\abs f}}{\grad f}.
\end{multlined}
\end{align*}
One has $\grad {\abs f}^{-1} = - f^{-2} \grad {\abs f}$ and $\extdiff\norm{\grad f}^2 = 2 \innerp*{\levici \grad f}{\grad f} = 2 \imath_{\grad f}\Hess(f)$, so $\grad \norm{\grad f}^2 = 2 \levici_{\grad f}{\grad f}$, so $\grad \norm{\grad f} = \levici_{\grad f}{\grad f} \norm{\grad f}^{-1}$, so
\begin{equation*}
\innerp*{\grad \norm{\grad f}}{\grad f} = \Hess(f)(\grad f, \grad f) \norm{\grad f}^{-1}.
\end{equation*}
Therefore,
\begin{multline*}
\div \left( F \circ \frac{\grad f}{f} \right) =
\left( G \circ \frac{\norm{\grad f}}{\abs f} \right) \left( \frac{\lapla f}{f} - \frac{\norm{\grad f}^2}{f^2} \right) +\\
\sigma_f
\left( G' \circ \frac{\norm{\grad f}}{\abs f} \right) \left( \frac{\Hess(f)(\grad f, \grad f)}{f^2 \norm{\grad f}} - \frac{\norm{\grad f}^3}{f^3} \right)
\end{multline*}
which yields the result.

\bigskip
\noindent
\textbf{Proof of 3.}
Let $g$ be as in the first part of the statement.
Apply (2.) to the function $G$ defined by $G(x) \coloneqq \frac{g(x)}{\sqrt{1+x^2}}$.
One has $G \circ \frac{\norm{\grad f}}{\abs f} = \left( g \circ \frac{\norm{\grad f}}{\abs f} \right) \frac{\abs f}{\eta_f}$ and
$G' \circ \frac{\norm{\grad f}}{\abs f} = \left( g' \circ \frac{\norm{\grad f}}{\abs f} \right) \frac{\abs f}{\eta_f} - \left( g \circ \frac{\norm{\grad f}}{\abs f} \right) \frac{f^2 \norm{\grad f}}{\eta_f^3}$.
Therefore,
\begin{multline*}
\div \left( F \circ \frac{\grad f}{f} \right) =
\sigma_f \left( g \circ \frac{\norm{\grad f}}{\abs f} \right)\\
\left( \frac{f}{\eta_f} \left( \frac{\lapla f}{f} - \frac{\norm{\grad f}^2}{f^2} \right) -
\frac{f^2 \norm{\grad f}}{\eta_f^3}
\left( \frac{\Hess(f)(\grad f, \grad f)}{f^2 \norm{\grad f}} - \frac{\norm{\grad f}^3}{f^3} \right)
\right) +\\
\left( g' \circ \frac{\norm{\grad f}}{\abs f} \right) \frac{1}{\eta_f}
\left( \frac{\Hess(f)(\grad f, \grad f)}{f \norm{\grad f}} - \frac{\norm{\grad f}^3}{f^2} \right)
\end{multline*}
which yields the result after simplification.
As for integrability, in local coordinates adapted to the boundary, one has $f(x) \sim_0 A/x_0$, so it is sufficient to prove that $g'(x^{-1}) x^{-2}$ is locally integrable at 0.
By change of variable, this amounts to the integrability of $g'$ at $+\infty$.

\bigskip
\noindent
Let $g$ be as in the second part of the statement.
Apply (2.) to the function $G$ defined by $G(x) \coloneqq \frac{g(x)}{x}$.
One has $G \circ \frac{\norm{\grad f}}{\abs f} = \left( g \circ \frac{\norm{\grad f}}{\abs f} \right) \frac{\abs f}{\norm{\grad f}}$ and
$G' \circ \frac{\norm{\grad f}}{\abs f} = \left( g' \circ \frac{\norm{\grad f}}{\abs f} \right) \frac{\abs f}{\norm{\grad f}} - \left( g \circ \frac{\norm{\grad f}}{\abs f} \right) \frac{f^2}{\norm{\grad f}^2}$.
Therefore,
\begin{multline*}
\div \left( F \circ \frac{\grad f}{f} \right) =
\sigma_f \left( g \circ \frac{\norm{\grad f}}{\abs f} \right)\\
\left(
\frac{f}{\norm{\grad f}} \left( \frac{\lapla f}{f} - \frac{\norm{\grad f}^2}{f^2} \right) -
\frac{f^2}{\norm{\grad f}^2}
\left( \frac{\Hess(f)(\grad f, \grad f)}{f^2 \norm{\grad f}} - \frac{\norm{\grad f}^3}{f^3} \right)
\right) +\\
\left( g' \circ \frac{\norm{\grad f}}{\abs f} \right)
\left( \frac{\Hess(f)(\grad f, \grad f)}{f \norm{\grad f}^2} - \frac{\norm{\grad f}^2}{f^2} \right)
\end{multline*}
which yields the result after simplification.
The integrability condition is the same as for the first case.

\bigskip
\noindent
\textbf{Proof of 4.}
The three formulas are obtained respectively by letting $g = 1$ in the first case of~(3) and letting $g = \frac2\pi \arctan$ and $g = \tanh$ in the second case of~(3).
\end{proof}

\section{Obtaining continuous integrands}

As remarked above, the only non-continuous terms of the integrands are of the form
\begin{equation}
\sigma_f h \left( \Hess(f)(\grad f, \grad f) - (\lapla f) \norm{\grad f}^2 \right)
\end{equation}
with $h \in \calC^1(M, \R)$, namely, $h = \left( g \circ \frac{\norm{\grad f}}{\abs f} \right) \eta_f^{-3}$ and $ h= \left( g \circ \frac{\norm{\grad f}}{\abs f} \right) \norm{\grad f}^{-3}$.
This is dealt with in~\cite{nodal} (in the case of Equation~\eqref{eq:algebraic} on the flat torus) using an integration by parts.
The same method extends to compact Riemannian manifolds as follows.
One has
\begin{equation*}
\Hess(f)(\grad f, \grad f) - (\lapla f) \norm{\grad f}^2 =
\innerp*{\grad f}{\levici_{\grad f} \grad f - (\lapla f)\grad f}
\end{equation*}
Therefore,
\begin{equation*}
\sigma_f h \left( \Hess(f)(\grad f, \grad f) - (\lapla f) \norm{\grad f}^2 \right) =
\innerp*{\grad{\abs f}}{h \left( \levici_{\grad f} \grad f - (\lapla f)\grad f \right)}.
\end{equation*}
We temporarily assume that $f$ is of class $\calC^3$ and we use the fact that $\div \left( \abs{f} h (\levici_{\grad f} \grad f - (\lapla f)\grad f) \right)$ has a vanishing integral on $M$ (by the standard divergence theorem).
Therefore,
\begin{multline*}
\int_M \sigma_f h \left( \Hess(f)(\grad f, \grad f) - (\lapla f) \norm{\grad f}^2 \right) \vol_M =\\
\int_M \abs{f} \div \left( h \left( (\lapla f)\grad f - \levici_{\grad f} \grad f \right) \right) \vol_M.
\end{multline*}
One has $\div((\lapla f)\grad f) = (\lapla f)^2 + \innerp{\grad \lapla f}{\grad f}$.
The Bochner formula yields
\begin{align*}
\div \left( \levici_{\grad f} \grad f \right)
&= \div \left ( \frac12 \grad \norm{\grad f}^2 \right) \\
&= \frac12 \lapla \norm{\grad f}^2\\
&= \innerp*{\grad \lapla f}{\grad f} + \norm{\Hess f}^2 + \Ric(\grad f, \grad f)
\end{align*}
where the norm of the Hessian is the Hilbert--Schmidt norm.
Therefore, the third derivatives cancel out.
Since $\calC^2(M, \R)$ is dense in $\calC^3(M, \R)$ for the $\calC^2$-topology (see for instance~\cite{hirsch}) and the involved quantities are continuous in this topology, one has, for any $f$ of class $\calC^2$,
\begin{multline}
\int_M \sigma_f h \left( \Hess(f)(\grad f, \grad f) - (\lapla f) \norm{\grad f}^2 \right) \vol_M =\\
\int_M \abs{f}
\left(
h \left(
(\lapla f)^2 - \norm{\Hess f}^2 - \Ric(\grad f, \grad f)
\right) +
\right.\\\left.
\innerp*{\grad h}{(\lapla f)\grad f - \levici_{\grad f} \grad f} \right) \vol_M.
\end{multline}

For example, one has $\grad{\eta_f^{-3}} = -3 \eta_f^{-5} \left( f \grad f + \levici_{\grad f} \grad f\right)$, so Formula~\eqref{eq:algebraic} becomes
\begin{multline}
\label{eq:lipschitz}
\vol(\calZ_f) = \frac12 \int_M
\left(
\frac{\abs f}{\eta_f^3} \left( \norm{\grad f}^2 - \abs f \lapla f + (\lapla f)^2 - \norm{\Hess f}^2 - \Ric(\grad f, \grad f) \right) + {}
\right.\\\left.
3 \eta_f^{-5} \left( \vphantom{(\lapla f)^2} f \Hess(f)(\grad f, \grad f) + \Hess(f)(\grad f, \levici_{\grad f} \grad f) - {}
\right.\right.\\\left.\left.
(\lapla f)^2 (f \norm{\grad f}^2 + \Hess(f)(\grad f, \grad f)) \right)
\vphantom{\frac{\abs f}{\eta_f^3}} \right)
\vol_M.
\end{multline}
This is the version of~\cite[Prop.~7]{nodal} for compact Riemannian manifolds.

These formulas can also be written in terms of the tracefree Hessian.
Recall that $\Hess^0(f) = \Hess(f) - \frac{\lapla f}{n} \id$.
One easily sees that a tracefree linear map is Hilbert--Schmidt-orthogonal to the identity, so $\norm{\Hess f}^2 = \frac{(\lapla f)^2}{n} + \norm{\Hess^0 f}^2$.

\end{document}